\documentclass[preprint,12pt]{elsarticle}
\pdfoutput=1
\usepackage{amsmath,amssymb,amsfonts}
\usepackage{algorithmic}
\usepackage{mathrsfs}
\usepackage{graphicx}
\usepackage{appendix}
\usepackage{textcomp}
\usepackage{float}
\usepackage{xcolor}
\usepackage{subfigure} 
\usepackage{diagbox}
\usepackage{hyperref}
\usepackage{multirow}
\usepackage{caption}
\usepackage{mathabx}
\usepackage{algorithm, algorithmic}

\usepackage{lineno}
\begin{document}
\begin{frontmatter}
\title{The Predictability of Stock Price: Empirical Study on Tick Data in Chinese Stock Market}

\author[a,b]{Yueshan Chen}
\ead{cys515@mail.ustc.edu.cn}
\author[b]{Xingyu Xu}
\ead{xingyuxu@mail.ustc.edu.cn}
\author[b]{Tian Lan}
\ead{lantian1@mail.ustc.edu.cn}
\author[b,c]{Sihai Zhang\corref{cor1}}
\ead{shzhang@ustc.edu.cn}
\address[a]{School of Cyberspace Science and Technology, University of Science and Technology of China} 
\address[b]{Key Laboratory of Wireless-Optical Communications, CAS}
\address[c]{School of Micro-Electronics, University of Science and Technology of China}

\cortext[cor1]{Corresponding author}

\begin{abstract}
Whether or not stocks are predictable has been a topic of concern for decades.
The efficient market hypothesis (EMH) says that it is difficult for investors to make extra profits by predicting stock prices, but this may not be true, especially for the Chinese stock market.
Therefore, we explore the predictability of the Chinese stock market based on tick data, a widely studied high-frequency data.
We obtain the predictability of 3, 834 Chinese stocks by adopting the concept of true entropy, which is calculated by Limpel-Ziv data compression method.
The Markov chain model and the diffusion kernel model are used to compare the upper bounds on predictability, and it is concluded that there is still a significant performance gap between the forecasting models used and the theoretical upper bounds.
Our work shows that more than 73\% of stocks have prediction accuracy greater than 70\% and RMSE less than 2 CNY under different quantification intervals with different models.
We further take Spearman's correlation to reveal that the average stock price and price volatility may have a negative impact on prediction accuracy, which may be helpful for stock investors.
\end{abstract}

\end{frontmatter}

\section{Introduction}
\label{sec:introduction}

The stock market has gradually become an important part of modern society, which has a fundamental impact on human life.
Stock ownership may help companies access capital more quickly and investors enjoy the benefits of corporate development, which implements more efficient resource allocation and ultimately drives faster social development.
But due to its high risk and high return characteristics, stock investments require effective analysis and reliable prediction which has extremely important theoretical significance and application value.
The efficient market hypothesis (EMH) assumes that the stock market has perfect law, good function, high transparency, full competition, where all valuable information is timely, accurate, fully reflected in the share price, including current and future value, unless there is market manipulation.
Thus investors can not obtain excess profits higher than the market's through the analysis of previous price \cite{1965Fama}.
Unfortunately, such strict conditions make semi-strong efficiency and strong efficiency almost impossible, even in well-capitalised markets \cite{grossman1980impossibility}.
For example, the American market is verified to be weakly efficient in both ordinary and special times \cite{sanchez2020testing}, while the Chinese stock market is not yet an efficient market \cite{2003Research}.
So this implies that current stock markets in the real world are somehow predictable.

Before the 1990s, most studies on financial time series were conducted on daily, weekly or monthly data, which is usually called low-frequency data in the field of financial econometrics \cite{2012Using}.
In recent years, with the development of computing tools and methods, the cost of data recording and storage has been greatly reduced, making it possible to study financial data with much higher frequency \cite{2017QuantCloud}.
In the financial markets, data collected in hours, minutes or seconds is high-frequency data, also known as intra-day data.

Tick data is one type of high-frequency data that is sampled in seconds, which is in fact a snapshot of stock market trading.
Stock exchanges send each snapshot to the market in real time, including current price, highest, lowest, trading volume, etc.\cite{nagy2012partitional}.

Up to now, there have been many kinds of prediction models for these types of stock price series data.
But the question of how predictable a certain stock market, or all stock markets, which motivates many researchers to spend numerous energy on this attracting issue.
The real entropy is adopted to measure the predictability of daily stock prices data, which is verified to be around 2.75, while the entropy of randomly shuffled variants of the original data achieve 2.90 \cite{2007On}. 
This validates that the stock time series is not completely random, but weakly time-dependent.
The stock price data with one minute granularity is further validated to have higher predictability than that with one day granularity \cite{fiedor2014frequency}.
Thus we raise two questions: (1) Does higher-frequency data have higher predictability? (2) What is the maximum prediction accuracy of this kind of data?

To explore these two issues, this paper adopts the concept of real entropy from the field of human location prediction into the tick financial data prediction.
Up to now, real entropy has been used to measure the uncertainty in many forecasting fields, such as wireless user mobility prediction \cite{liu2019diffusion}, wireless traffic prediction \cite{2021Can}, taxi demand prediction \cite{2016taxi}, electronic health records \cite{dahlem2015predictability}, cyberattacks \cite{chen2015spatiotemporal} and human-interest dynamics \cite{2013Emergence}.
All these studies help relevant industries to solve the theoretical problem of time series predictability.
Thus we use this concept to investigate the predictability in the field of financial stock price prediction based on Chinese stock data.

The contributions of this work are as follows:
\begin{itemize}
\item We use real entropy to calculate the predictability of tick data.
By using different quantification interval settings, the state spaces of stock time series are established based on which the entropy along with predictability are obtained and analysed. 
Our result based on the three-second tick data including 3,834 stocks with 230 million records shows that 74\% of stocks have a real entropy of less than 2, which means most stocks in Shanghai and Shenzhen markets have high predictability.

\item We evaluate the Markov Chain model and Diffusion Kernel model based on tick data and demonstrate that the accuracy of both algorithms still have significant gaps of 0.090 and 0.126 to the upper bound, respectively.
We also provide the trade-off of price quantification in $T = 0.01$ and $T = 0.05$ using predictability which helps evaluate the prediction precision and accuracy of stock price.
 
\item We find that the predictability and prediction accuracy of all stock prices shows severe inherent differences.
We discuss the six features related to each stock to do a correlation analysis with accuracy, including the industry of the company to which the stock belongs, region, company scale, and lifetime of the stock, as well as the average price and historical volatility of the stock price series.
The correlation coefficients of stock price and volatility with accuracy are -0.72 and -0.58, which are closely and negatively correlated.
\end{itemize}

The remainder of this paper is organized as follows. 
Sec.\ref{sec:related work} summarizes the stock price forecasting models in recent years.
Sec.\ref{sec:concepts} briefly introduces the fundamental concepts and methods used in this research. 
In Sec.\ref{sec:dataset}, we describe the data set used in this experiment. 
Then, we use two models for prediction, and the results and analyses are presented in Sec.\ref{sec:performance}. 
Besides, the experimental results are further discussed and we make supplementary experiments and analyses in Sec.\ref{sec:discussions}. 
Finally, the conclusion is drawn in Sec.\ref{sec:conclusions}.

\section{Related work}
\label{sec:related work}

In recent years, there has been an increasing amount of literature on stock price forecasting.
According to these modeling theories, these prediction models can be divided into two categories, one is time series analysis based on statistical principles, and the other is the new artificial intelligence prediction method \cite{2021Artificial}.

The time series model mainly uses statistical software to extract the favorable information of historical stock price to construct the prediction model, such as generalized autoregressive conditional heteroskedasticity (GARCH) and Autoregressive Integrated Moving Average model (ARIMA).
GARCH is widely used for time series prediction \cite{2019Literature}, while ARIMA is proved to be particularly effective in the cases of short-term forecasts \cite{2020Forcasting}.
Such models are based on the assumption that the investigated financial time series are generated by a linear process and model the time series in order to predict the future value of the series.
However, stock time series data are inherently highly noisy, nonlinear, complex, dynamic, nonparametric and chaotic \cite{2021survey}.
As such, it is not suitable for traditional statistical techniques to model the complexity and non-stationarity of stock markets, so that many researchers use nonlinear models of artificial intelligence (AI) or combine statistical methods with these models.


These AI models include Support Vector Machine (SVM), Genetic Algorithm (GA), Convolutional Neural Network (CNN), Artificial Neural Network (ANN), Long Short Term Memory Network(LSTM), Fuzzy Logic (FL) etc.
Among them, ANN is one of the most frequently used models to deal with nonlinear problems, an ANN model combined with a nonlinear autoregressive model has been shown to have better predictive performance than Exponential GARCH (E-GARCH, a competing asymmetric conditional variance model with the superior predictive performance), but the training of neural network is time-consuming and easy to overfitting \cite {2021ANN}.
CNN is also widely used as a strong benchmark for any innovative machine learning model, but LSTM may be a better choice for prediction accuracy \cite{nabipour2020deep}.
This conclusion is also confirmed by work called LSTMLI and SSACNN that uses options and futures along with historical stock data as input to LSTM and CNN for prediction \cite{SSACNN}.
And the work concluded that LSTMLI has higher accuracy than SSACNN and other models like SVM \cite{LSTMLI}.
However, even with the improvement of LSTM framework, one single model can obtain the accuracy of predicting the rise and fall of low-frequency data in the stock market to only 55\%-60\%.
Therefore, most research on stock price prediction models are aimed at the improvement of LSTM or the mixture of various machine learning models, and all have achieved better results than single model \cite{2020Financial, 2018Which}.

However, multi-model hybrid prediction methods require a large number of hyperparameters and considerable effort to optimize, while Automated Machine Learning (autoML) can solve this problem by automatically finding the suitable machine learning model and optimizing it \cite{autoML}. 
There has been some research on time series forecasting problems on autoML systems such as GluonTS \cite{GluonTS}, AutoAI-TS \cite{AutoAI-TS} etc., which enables multiple AI models to automate the time series forecasting process and can achieve better performance than deep learning models on most non-financial datasets.
However, due to the drawback of limited search space in the financial data processing pipeline \cite{GECCO}, autoML does not perform well on datasets in the financial domain, 
e.g., an AutoML approach named BOHB can obtain higher precision than the manual DL method on the voltage dataset and the methane dataset, but not as good on the New York Stock Exchange (NYSE) and Johannesburg Stock Exchange (JSE) datasets \cite{HOBO}.
AutoML for time series is still in the evolution stage and requires the efforts by researchers to model more adaptable AutoML frameworks that can accommodate datasets in different domains and conduct empirical studies.

Despite the variety of methods available today for forecasting stock data, there is still a lack of theoretical guidance for quantitative analysis of prediction model performance bounds in this field.
Therefore, it is necessary to analyse the time series of stocks and calculate the maximum value that can theoretically achievable in terms of forecasting performance, which can be used to guide the improvement of advanced forecasting models.

\section{Concepts and Methods}
\label{sec:concepts}

In this section, the fundamental concepts and methods used in this research are briefly introduced.

\subsection {Real Entropy}

Entropy generally refers to the degree to measure the uncertainty of one random variable.
But neither the Shannon entropy nor random entropy can capture the characteristics of the time series representing the moving positions.
However, the real entropy measurement not only considers the occurrence frequency of states, but also digs the occurrence order of states and the time of each state stay, so that the complete spatial and temporal information in the complete time series are excavated \cite{song2010limits}.
For easy reading, the definition of real entropy is briefly introduced as follows, and for more details, please refer to \cite{song2010limits}.

Suppose a historical sequence $T=\{X_1, X_2, \dots, X_n\}$, then the real entropy $S$ is defined as
\begin{align} 
\label{S}
S = -\sum_{T^\prime\subset{T}} p(T^\prime)\log_2(T^\prime)
\end{align}
where $p(T^\prime)$ represents the probability of finding a subsequence $T^\prime$ in the trajectory $T$.
For the time series with length $n$, the real entropy can be estimated by Lempel-Ziv data compression method as follows
\begin{align} 
\label{S^{est}}
S^{est} = (\frac{1}{n}\sum_{i}^{n}\Lambda_i)^{-1}\ln n
\end{align}
where $\Lambda_i$ is the shortest length $k$ such that the subsequence starting from position $i$ with length $k$ does not appear previously as a continuous subsequence of $ \{X_1, X_2, \dots, X_{i-1}\}$.

\subsection {Predictability}
\label{subsec:pred}

Predictability indicates the probability that one optimal prediction algorithm (theoretically) predicts the user’s future state correctly \cite{2021Can}.
If the predictability can be calculated, the gap between the existing prediction methods and the theoretical optimal algorithm can be recognized.
From the estimation of the real entropy by Eq.\ref{S^{est}}, the upper bound of predictability $\varPi$ can be obtained using the Fano inequality, as follows: \cite{song2010limits}

\begin{equation}
\label{Fano_S}
S=H(\varPi)+(1-\varPi)\log_2(N-1)
\end{equation} 

\begin{align}\label{H(Pi^max)}
H(\varPi) = -\varPi\log_2{(\varPi)}-(1-\varPi)\log_2{(1-\varPi)}
\end{align} 
where $N$ denotes the number of distinct states appeared in sequence $T$. 

\subsection {Prediction Models}
\label{subsec:models}

In this paper, two models, second-order Markov Chain model and second-order Diffusion Kernel model, are used to predict the stock price.

The Markov Chain (MC) model is a stochastic process in which the historical state has nothing to do with the future state under the given current information \cite{gagniuc2017markov}. 
The $k$-order Markov Chain model has higher complexity and lower accuracy when $k>2$, so first-order and second-order MC are usually selected to realize the location prediction tasks.
The second-order Markov model means that the state of the model at time $t$ is only related to the state at time $t-1$ and time $t-2$, and independent of other states.

The Diffusion Kernel (DK) model is similar to the neural network prediction model and gradient descent operations are involved in the calculation without requiring large amounts of training data. 
The concept of DK model is to map the mobility trajectory into diffusion process in a continuous space which embeds the visited locations.
The model is illustrated in Algorithm.\ref{alg}, and for more details, please refer to \cite{liu2019diffusion}.

\begin{algorithm}[htb]
	\renewcommand{\algorithmicrequire}{\textbf{Require:}}
	\caption{Prediction based on Diffusion Kernel model}
	\label{alg}
	\begin{algorithmic}[1]
		\REQUIRE a state space $L$, the training set of a sequence set of state $Traj_l$ with time stamp $t$, the set of coordinate of nodes in an Euclidean space $Z$ \\
		\STATE $\forall l\in L,z^{\tau}_{l} \gets random,\tau \gets 0$\\
		\WHILE{$\tau < T$}
		\FORALL{$traj\in Traj_l, l_i\in traj$, $l_j\in L$ with $t(l_i)<t(l_j)$ or $l_j \in \widebar{traj}$}
		\STATE $d_i \gets \|z_{s^c}-z_{l_i}\|^2$\\ 
		\STATE $d_j \gets \|z_{s^c}-z_{l_j}\|^2$\\ 
		\STATE $\delta_d \gets d_j-d_i$\\ 
		\IF{$\delta_d<1$}
		\STATE $z_{l_i}^{(\tau +1)}\gets z_{l_i}^{(\tau)}+\alpha(\tau)\times2(z_{s^{traj}}^{(\tau)}-z_{l_i}^{(\tau)})$\\
		\STATE $z_{l_j}^{(\tau +1)}\gets z_{l_j}^{(\tau)}+\alpha(\tau)\times2(z_{s^{traj}}^{(\tau)}-z_{l_j}^{(\tau)})$\\
		\STATE $z_{s^{traj}}^{(\tau +1)}\gets z_{l_i}^{(\tau)}+\alpha(\tau)\times2(z_{l_i}^{(\tau)}-z_{l_j}^{(\tau)})$\\
		\ENDIF
		\ENDFOR
		\STATE $\tau \gets \tau +1$
		\ENDWHILE
	\RETURN $Z$
	\end{algorithmic}  
\end{algorithm}

Generally speaking, the prediction performance of the DK model increases with the higher orders, but with more computation complexity and costs.
In this work, in order to ensure the timeliness of prediction results of 3-seconds tick data, the second-order DK model is adopted.

\section {Data Sets}
\label{sec:dataset}

\begin{figure}[!h]
	\centering
	\includegraphics[scale=1.2]{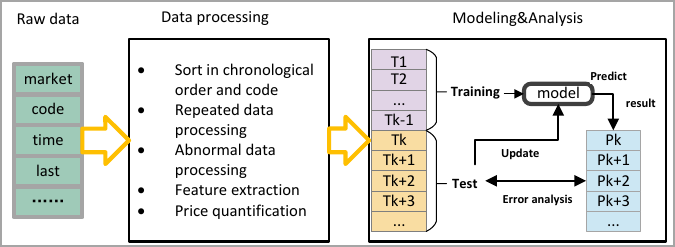}
	\caption{The architecture workflow of data processing and modeling}
	\label{fig:workflow}
\end{figure}

This section describes the dataset used in this work and the pre-processing steps for the dataset.
The architecture workflow is shown in Fig.\ref{fig:workflow}.

We use two datasets in this work, and they are the tick data of Shanghai Stock Market and Shenzhen Stock Market in January 2021, and Shenzhen Stock Market in July 2022, which are called dataset A and B, respectively.
Tick Data is the buy and sell orders that the exchange sends to the stock market in the active Order Book for each stock.
The domestic exchange takes a snapshot in each three seconds, which is actually like a snapshot of each stock with current price, highest, lowest, trading volume, transaction amount and other market features in the latest three seconds.
Daily continuous bidding periods in the two stock markets include two hours in the morning and two hours in the afternoon, which leads the number of snapshots to be about 3,800.
There is more than 2 GB of market-wide snapshot data about the whole stock per day.

The dataset A has a total of 4,147 stocks, and each stock is gathered in a time granularity of three seconds. 
Among these stocks, 1,800 stocks are from the Shanghai stock market, while the rest are from the Shenzhen stock market.
The Dataset A includes 18 trading days, and each record includes 56 features such as stock code, time, open price, close price and last-price, etc. 
The last price of the stock means the newest real-time price of the stock in the process of real-time change.
The stock price fluctuates ceaselessly, so the last price is also dynamic, which can reflect the stock change situation.
Thus, the last-price of each record is extracted as the last-price time series of all stocks for the following analysis.

The dataset B can be downloaded in \url{https://www.kaggle.com/datasets/chenyueshan/shenzhen-tick-data}.
We pick 2,324 stocks in the validation dataset that have the same code as the dataset A in the Shenzhen stock market, and construct a time series for each stock.
Due to the difference in data collection platforms, the validation dataset was not collected as frequently as the main dataset, at around 10-second intervals.
This dataset has 21 trading days with 27 features, including time, turnover, volume, etc.
In this dataset we chose turnover as the target for prediction.

The prediction target is the last-price and turnover at the next moment.
Since tick data takes seconds as time granularity, we can approximately regard the prices as continuously changing prices.
Consequently, an quantification interval should be adopted to quantify the stock price, making the stock price within a range into a state.
We can set different quantification intervals according to different predicted demands.
Since the precision of stock price in the data set is 0.01 CNY, the minimum quantification interval we can set is 0.01 CNY, that is, not to quantify the stock sequence.
We initially set two quantification intervals, 0.01 CNY and 0.05 CNY, respectively.

\begin{figure}
    \hspace{-3mm}
	\subfigure[number of valid states]{
	    \centering
	    \label{fig:compare_statenum}
		\includegraphics[width=0.55\linewidth]{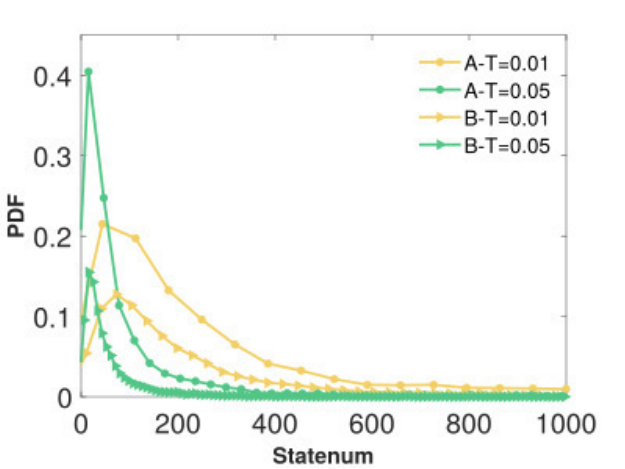}}
	\hspace{-10mm}
	\subfigure[real entropy]{
	    \centering
	    \label{fig:compare_entropy}
		\includegraphics[width=0.55\linewidth]{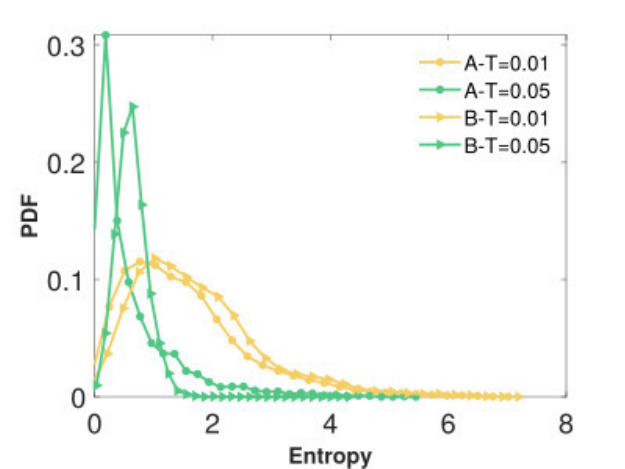}}
	\caption{The subfigure\ref{fig:compare_statenum} and\ref{fig:compare_entropy} show the distribution of valid state space and the distribution of real entropy of stocks in $T=0.01$ and $T=0.05$ from dataset A and B, respectively.}
	\label{dataset} 
\end{figure}

We then take two preprocess operations on the quantifized data, filtering stocks whose sequence length is too short and state space is less than 10.
Because these stocks have missing data or they have been delisted, there is no actual predictive significance.
And then we make predictions for the remaining stocks. 
Fig.\ref{fig:compare_statenum} shows the distribution of the number of valid states for the complete sequence of these stocks, including only stocks with a state space size of [0,1000].
Here, a valid state is defined as the state in which the times of such a state occur at least once.
It indicates that the distribution of the number of valid states and real entropy of stocks calculated from the two data sets is similar, so as to avoid verbosity, we will mainly analyse dataset A.
Obviously, stocks have smaller state space when they are quantified by larger intervals.
Although the X-axis ranges from 0 to 1,000, in fact, there are also 209 stocks with more than 1,000 states at 0.01 CNY quantification intervals, while 58 stocks at 0.05 in dataset A.
Such stocks are relatively few and their number of valid states is discretely distributed. 
In order to clearly show the number of state space of most stocks, we do not show these stocks in Fig.\ref{fig:compare_statenum}, but the data set used in the experiment still includes these stocks.

While modeling, we set the first day as a training set, and the remaining days as a test set.
After training the model with the training set, we update the model after each prediction during the testing phase, so changing the ratio between the training set and the test set has no influence on the prediction result.

We use the Limpel-Ziv data compression method to calculate the real entropy of each stock, and the real entropy distribution of the last price is exhibited in Fig.\ref{fig:compare_entropy}.
In dataset A, about 74\% of stocks ($T=0.01$) and 87\% stock ($T=0.05$) have very low real entropy less than 2, respectively, which indicates that based on the existing historical price series, the next three seconds price of most stocks can be found in fewer than $2^2$ states.
This observation implies that the uncertainty of the forecast on the tick data is quite low, thus the stock price may be easy to predict.

\section {Performance Evaluation}
\label{sec:performance}

In this section, we first introduce the performance metric.
Then we use the MC and DK models to predict the stock price with different intervals, observe the gap between the predictability upper bound and the prediction accuracy of models.
We also present the prediction error measured by RMSE and analyzed the relationship between prediction accuracy and RMSE.
Finally, considering the impact of stock price, we perform a supplementary experiment to divide the quantification interval according to the high and low price of the stock price series, then evaluate the results by accuracy and the ratio of DK-RMSE to stock price.

\subsection{Performance Metric} 

As for performance metrics, we use accuracy $ACC$ and Root Mean Squared Error (RMSE), which are introduced as follows:

\begin{align} 
	\label{ACC}
	ACC = \frac{{c_{test}}}{{n_{test}}}
\end{align}
\vspace{-5mm}
\begin{align}
	RMSE = \sqrt{\frac{1}{n}\sum_{t=1}^{n} (y_{t}-y_{t}^{\prime})^2}
	\label{RMSE}		
\end{align}

where $c_{test}$ is the number of states predicted correctly in test set, $n_{test}$ is the length of test set sequence, $y_{t}$ is the actual stock price at time $t$, and $y_{t}^{\prime}$ is the predicted price at time $t$.

\subsection{Accuracy} 

\begin{figure}
    \hspace{-3mm}
	\subfigure[Dataset A]{
	    \centering
		\includegraphics[width=0.55\linewidth]{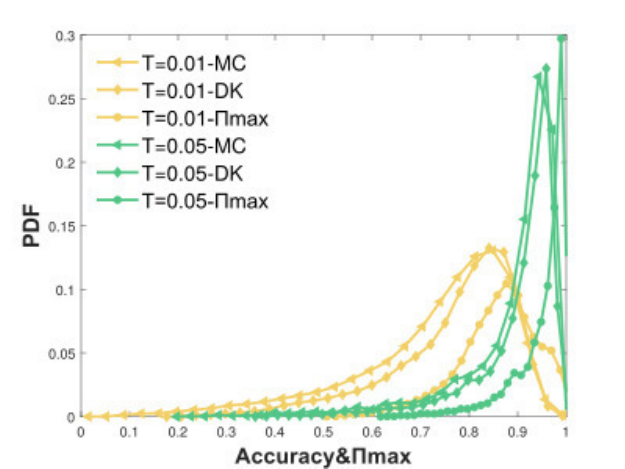}}
    \hspace{-10mm}
	\subfigure[Dataset B]{
	    \centering
		\includegraphics[width=0.55\linewidth]{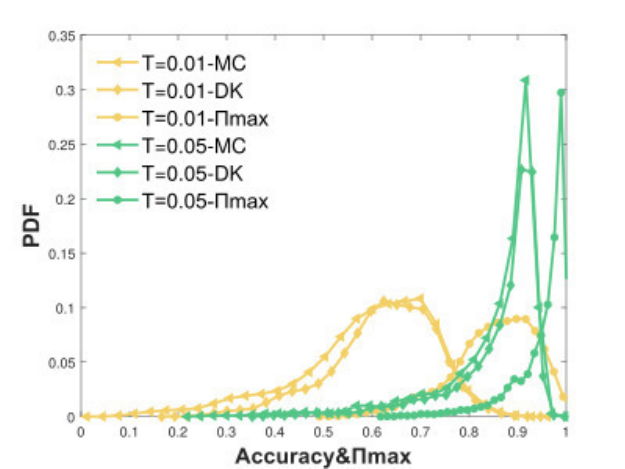}}
	\caption{Accuracy and $\Pi_{max}$ distribution of stocks in $T=0.01$ and $T=0.05$}
	\label{fig:compare_acc}
\end{figure}

Once we have the real entropy of each stock, we can obtain upper boundaries on predictability.
At the same time, we use the two models to forecast the tick data's last price and observe the gap between the two models’ prediction accuracy and the upper bound of predictability in this section.
Fig.\ref{fig:compare_acc} and Table.\ref{TABLE:ACC} demonstrates the ACC distribution results based on MC and DK models and the maximum predictability under $T = 0.01$ and $T = 0.05$.

\begin{table}[htb]
\setlength{\abovecaptionskip}{0mm}
\setlength{\belowcaptionskip}{2mm}
\caption{Arithmetic mean ACC and $\Pi_{max}$ in different intervals of dataset A and B}
\centering
\begin{tabular}{|c|c|c|c|c|c|c|}
\hline
		\multirow{2}*{Dataset}&\multicolumn{3}{|c|}{T=0.01} &\multicolumn{3}{|c|}{T=0.05} \\
		\cline{2-7}
		~  & MC &DK & $\Pi_{max}$& MC &DK & $\Pi_{max}$\\
\hline
A & 0.732 & 0.768 & 0.858 & 0.888 & 0.899 & 0.951 \\
\hline
B & 0.59 & 0.622 & 0.851 & 0.850 & 0.858 & 0.945\\
\hline
\end{tabular}
\label{TABLE:ACC}
\end{table}

The most important finding is that neither MC nor DK model achieves the theoretical upper bound of ACC, $\Pi_{max}$ when $T = 0.01$ and $T = 0.05$.
To be specific, when $T = 0.01$ the arithmetic mean ACC of all stocks are 0.732 and 0.768 with DK and MC models respectively, while $\Pi_{max}$ is 0.858 in Dataset A.
Note that, in this work only MC and DK models are verified, but according to the existing reports, other prediction models, such as BoF \cite{passalis2020temporal}, LSTM \cite{fischer2018deep}, are also not able to achieve such high accuracy.
So in summary, our finding indicates that current prediction performance of existing models on tick data in the Chinese market still have room to increase.

Then the ACC of MC and DK are compared, which indicates that DK's ACC is slightly higher than MC's both when $T = 0.01$ and $T = 0.05$.
That is because the DK model uses latent representation rather than the transition matrix, thus overcomes the limitation of the MC caused by the high complexity transition matrix in this case of a large state space \cite{liu2019diffusion}.

We also note here that, when $T = 0.05$, the ACC and prediction upper bound are higher than that when $T = 0.01$.
The reason for this is quite straightforward, because the size of state space is smaller when $T = 0.05$ and the entropy will be lower, therefore accuracy and $\Pi_{max}$ will be higher.

\subsection{RMSE}

\begin{figure}
    \hspace{-3mm}
	\subfigure[Dataset A]{
	    \centering
		\includegraphics[width=0.55\linewidth]{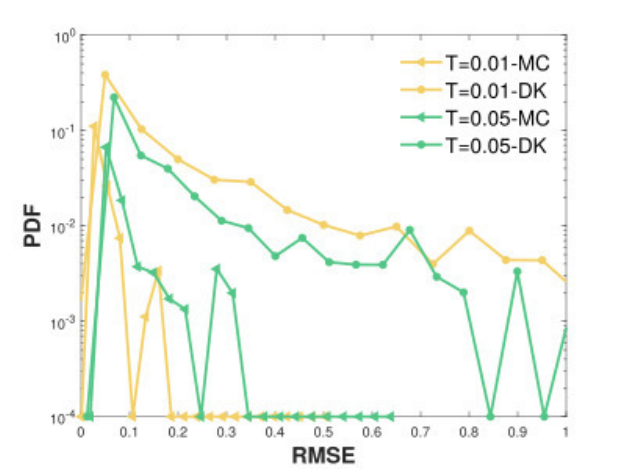}}
    \hspace{-10mm}
	\subfigure[Dataset B]{
	    \centering
		\includegraphics[width=0.55\linewidth]{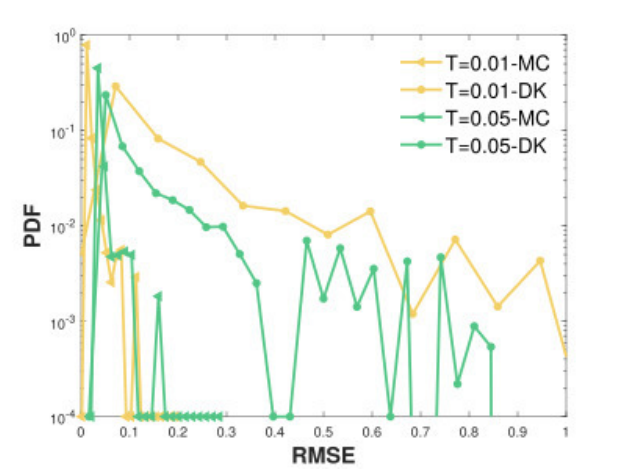}}
	\caption{RMSE distribution of stocks in $T=0.01$ and $T=0.05$ intervals in Dataset A and B. We intercepted the part with RMSE less than 1 CNY. In order to clearly compare the differences in different intervals and different models, we intercepted the part of RMSE less than 1 yuan here. For the complete RMSE distribution in dataset A, please refer to Fig.\ref{RMSE-acc}.}
    \label{fig:compare_rmse}
\end{figure}

RMSE is the popular metric to evaluate the performance of stock price prediction models, so we explore the RMSE performance under the influence of $T = 0.01$ and $T = 0.05$ on precision, which is shown in Fig.\ref{fig:compare_rmse} and Table.\ref{table:RMSE}.
Clearly, the RMSE is smaller in $T = 0.01$ than in $T = 0.01$ for both MC and DK models, which is easy to understand, because increasing the quantification interval can improve accuracy, although it leads to a decrease in precision.

\begin{table}
\setlength{\abovecaptionskip}{0mm}
\setlength{\belowcaptionskip}{2mm}
\caption{Arithmetic mean RMSE in different intervals of dataset A and B}
\centering
\begin{tabular}{|c|c|c|c|c|}
\hline
		\multirow{2}*{Dataset}&\multicolumn{2}{|c|}{T=0.01} &\multicolumn{2}{|c|}{T=0.05} \\
		\cline{2-5}
		~  & MC &DK & MC &DK \\
\hline
A & 0.015 & 0.312 & 0.035 & 0.211 \\
\hline
B & 0.013 & 0.149 & 0.033 & 0.100\\
\hline
\end{tabular}
\label{table:RMSE}
\end{table}

However, the results of the comparison between the two models are the opposite of what we expected.
As can be seen from Fig.\ref{fig:compare_acc}, DK's result is more accurate than MC's.
So theoretically its RMSE should be smaller than MC, but the actual results shown in Fig.\ref{fig:compare_rmse} are the opposite.
Here we focus our analysis on dataset A.
Just as it shows, the RMSE of all stocks in both quantification intervals is less than 0.64 when using the MC model, while the DK model performs significantly worse than MC.
Moreover, the statistically significant number of stocks with RMSE greater than 1 in the DK model is 215 and 134 at $T = 0.01$ and $T = 0.05$, respectively.
But statistically, there are 215 and 134 stocks whose RMSE is bigger than 1 with the DK model at $T = 0.01$ and $T = 0.05$, respectively.
Among these, two stocks, IMEIK and Kweichow Moutai, have quite big RMSEs, reaching 57 and 89 CNY respectively at $T = 0.01$.
Upon analysis, we noticed that these two stocks have a high stock price in common of 700 CNY and 2,000 CNY.
We then checked other stocks with big RMSEs and found that they are also traded at relatively high prices.
Therefore, we think that the prediction performance of the DK model is related to stock price, and it is verified in section\ref{sec:discussions}.
This explains why DK is more accurate than MC, but the RMSE is bigger, because DK is affected by the stock price but MC is not.

And when we adjust the quantification to $T = 0.05$, the RMSE of the high-priced stocks all decrease, such as IMEIK and Kweichow Moutai's RMSE decrease to about 41 and 66 CNY.
Therefore, when using the DK model to predict high-priced stocks, we can reduce the RMSE by using a larger quantification interval.

From the above analysis, we conclude that the RMSE of the DK model is related to the quantification interval and the stock price.
This is why the difference in MC results between dataset A and B is small, but the difference in DK model results is large. 
Because the stocks included in dataset B are all underpriced, with only 29 stocks having an average price of over 100 CNY, the RMSE for dataset B is smaller than that for dataset A.
Hence, when using the DK model, the appropriate quantification interval should be selected based on the stock price.

\subsection{Correlation of Accuracy and RMSE} 

\begin{figure}
    \hspace{-3mm}
	\subfigure[Dataset A]{
	    \centering
		\label{R-A-MC}
		\includegraphics[width=0.55\linewidth]{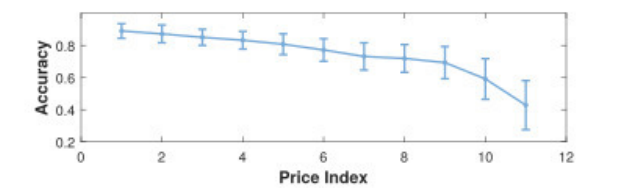}}
    \hspace{-10mm}
	\subfigure[Dataset B]{
	    \centering
		\label{R-A-DK}
		\includegraphics[width=0.55\linewidth]{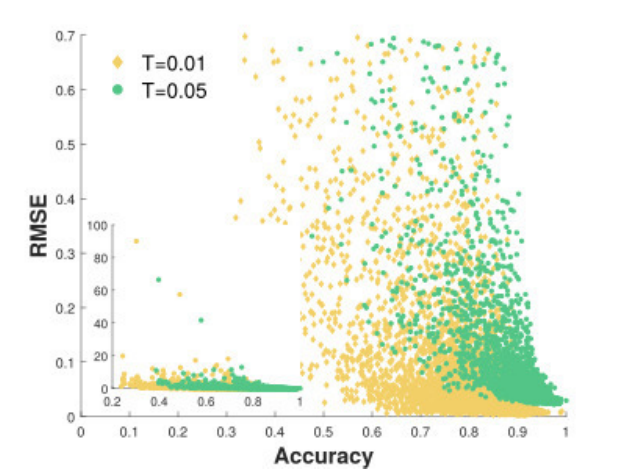} }
	\caption{The relationship between RMSE and accuracy with MC and DK model.}
	\label{RMSE-acc} 
\end{figure}

Fig.\ref{RMSE-acc} shows the correlation between RMSE and accuracy with MC model and DK model in Dataset A.
First of all, by comparing the two quantification intervals, we can still get the previous conclusion that setting a larger quantification interval leads to higher prediction accuracy but lower precision.
Likewise, we compare the two models and find that MC's RMSE is generally lower than that of DK's model, with DK's RMSE being related to the stock price.

However, combining Fig.\ref{R-A-MC} and the subplot of Fig.\ref{R-A-DK} to observe the overall distribution, we find that the general distribution trend is similar and can be divided into three parts, which are high accuracy with small RMSE, low accuracy with high RMSE, and low accuracy with low RMSE.
It is quite logical that the first two cases appear, and the reasons for them are obvious.
Then, we focus on why there are stock prediction results with low accuracy but small RMSE.
In the case of quantifying stocks, we consider a certain range of prices as a state.
And when we calculate the accuracy, we only determine whether the predicted state is consistent with the actual state, while the RMSE is calculated considering the distance between states.
If the predicted state is near the actual state in each prediction, low accuracy but high precision will occur.

Based on the distribution of the samples in Fig.\ref{RMSE-acc}, we believe that there is a relationship between accuracy and precision, but it requires complex mathematical modeling and theoretical analysis, which has not been covered in the current work. 
If we can get the relationship between them, we can combine the theory of calculating predictability by real entropy to get the upper bound of prediction precision of the prediction problem, that is, the lower limit of RMSE, and this is our prospect for future work.

\subsection {Supplementary Experiment}

Since we summarized in the previous section that the RMSE of the DK model is related to the stock price and quantification interval, we conducted a supplementary experiment with the DK model, using the dataset A which has more stock sample. 
Considering the range of the stock price sequence of the training set, the size of the state space is fixed for each stock to obtain the appropriate quantification interval.

In order to investigate what the appropriate size of the state space is, we set the size of the state space to 50, 100, and 150 respectively to conduct the experiment.
And we find that if the state space is set as 150, with the exception of several stocks, almost all stocks are quantificated by $T = 0.01$. 
And when their state space is set as 50, their quantification intervals are greater than $T = 0.05$.
As a result, we finally choose a state space of 100 ($SP=100$ for short) as the appropriate size and compare it with the $T = 0.01$ and $T = 0.05$.

\begin{figure}
    \hspace{-3mm}
        \subfigure[Accuracy and $\Pi_{max}$ distribution ]{
	    \label{fig:compare_DK_ACC}
		\includegraphics[width=0.55\linewidth]{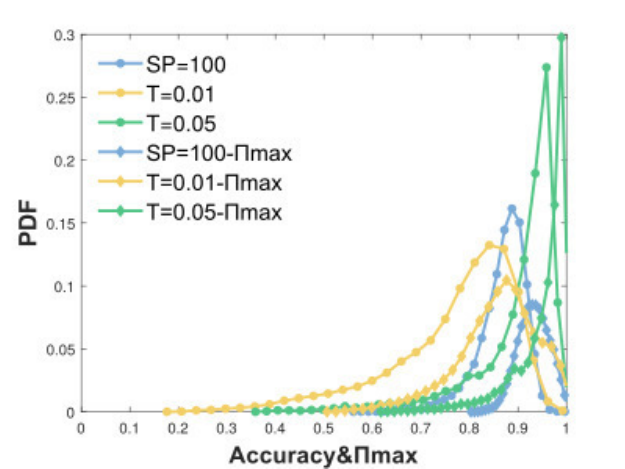}}
    \hspace{-10mm}
	    \subfigure[The ratio of DK-RMSE to price distribution]{
	    \label{fig:compare_ratio}
		\includegraphics[width=0.55\linewidth]{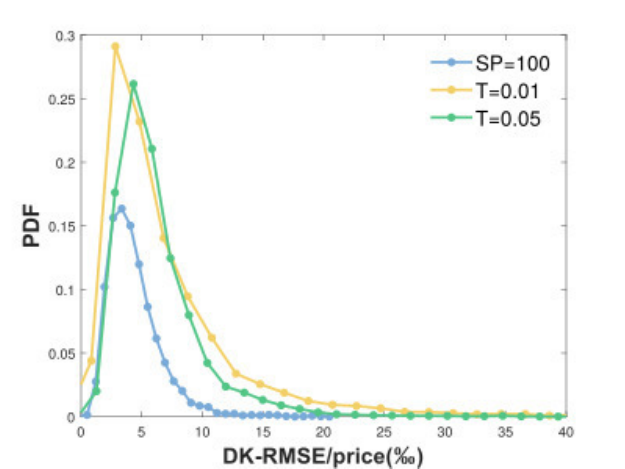} }
	\caption{The results of supplementary experiments in $T=0.01$, $T=0.05$ and $SP=100$. $SP = 100$ means the size of state space is fixed as 100.}
    \label{supplements} 
\end{figure}

As expected, both the accuracy curve and the $\Pi_{max}$ curve for the moderate state space lie between the other two cases as shown in Fig.\ref{fig:compare_DK_ACC}.
In particular, we find that the prediction accuracy of stocks at $SP = 100$ is distributed between 0.58 and 1, which is unlike the other two cases where there are stocks with accuracy less than 0.5.
Similarly, the upper bound of predictability for stocks is above 0.8, while many stocks have upper bounds in the range of [0.5, 0.8] when $T = 0.01$ and $T = 0.05$.
Combined with the above two findings, we can know that fixing state space can improve the accuracy and predictability of stocks with low prediction accuracy.

As for the RMSE results under three different quantification intervals, we consider the case where the same error has different importance for different levels of stock prices.
For example, an RMSE of 0.1 CNY has a deep impact on stocks with an average stock price of 1 CNY, but this error is insignificant for stocks with an average price of 1000.
To eliminate this difference, we set up an evaluation metric that uses the ratio of the predicted RMSE of the DK model to the average price of the stock price series, and re-observe the performance of DK in Fig.\ref{fig:compare_ratio}.
The arithmetic mean result of ratios for $SP = 100$, $T = 0.01$, and $T = 0.05$ are 4.35\textperthousand, 7.45\textperthousand, 6.58\textperthousand, respectively.
$SP = 100$'s result is obviously better than the other two cases.

\section{Discussions}
\label{sec:discussions}

The above analysis raises one key question: why do some stocks have higher predictability and prediction accuracy while others do not?
So we try to discuss and find certain characteristics related to stocks to indicate whether one stock is easy to predict in this section.

We arbitrarily choose six features, $category$, $region$, $scale$, $life$, $avgprice$ and $volatility$, related to the stock data and its belonging enterprises, four of which are profile of the enterprises while other two of which are price measurements. 

$category$ represents the industry to which the enterprise belongs.
For $category$, we classify 4,147 enterprises into twenty categories and number them according to Fortune China's new classification standards for Chinese industries.
The detailed information of these industries is shown in appendix.
$region$ means the geographical region where the enterprise is located, while 32 regions were numbered by Chinese provinces in appendix.
$scale$ denotes the number of employees in the enterprise.
$life$ represents how many years the enterprise has been traded on the stock market.
In our data set, the minimum life is one year and the maximum is 31 years.
Please refer to the appendix for more details.

$avgprice$ is the arithmetic mean last-price of the complete historical series and $volatility$ is a measure of the intensity of price movements.
The historical volatility $\sigma$ is calculated based on the definition\cite{andersen1998answering} as follows.

\begin{equation}\label{volatility}
\sigma = \sqrt{\frac{\sum(x_i-\bar{x})^{2}}{N-1}}
\end{equation} 

where $x_i=\ln(\frac{P_t}{P_{t-1}})$ is log-return, $\bar{x}=\frac{1}{N}\sum{x_i}$ is the mean of the return, $N$ means the length of historical price sequence $\{P_1, P_2, \dots, P_n\}$, and $P_t$ is the price at time $t$.

\begin{table}
\setlength{\abovecaptionskip}{0mm}
\setlength{\belowcaptionskip}{1mm}
\caption{ANOVA for region and industry}
\centering
\begin{tabular}{|c|c|c|c|c|c|c|}
\hline
Feature & F & p & SSB & SST &  ${{{\eta}^2}_p}$ \\
\hline
category & 2.179 & 0.002 & 0.676 & 67.768 & 0.010  \\
\hline
region & 4.267 & 0.000 & 0.818 & 34.083 & 0.024 \\
\hline
\end{tabular}
\label{ANOVA}
\end{table}

Table.\ref{ANOVA} shows the ANOVA (Analysis of Variance) results to investigate the effect of regions and categories on stock prediction accuracy.
Here, $F$ value is the statistic of F test, $p$ is the indicator of the significance of differences, $SSB$ represents the sum of squares between groups, and $SST$ is the sum of squares for total, ${{{\eta}^2}_p} = SSB / SST$ is the effect size partial eta square, which represents the magnitude of difference between groups.
The results may come to the conclusion that, in Chinese Shanghai and Shenzhen stock exchanges, the prediction accuracy of stock prices may be irrelevant to categories or regions. 
In other words, the industry and the location of the enterprises might not be the proper factor when predicting the prices.

\begin{figure}
	\centering
	\includegraphics[scale=0.6]{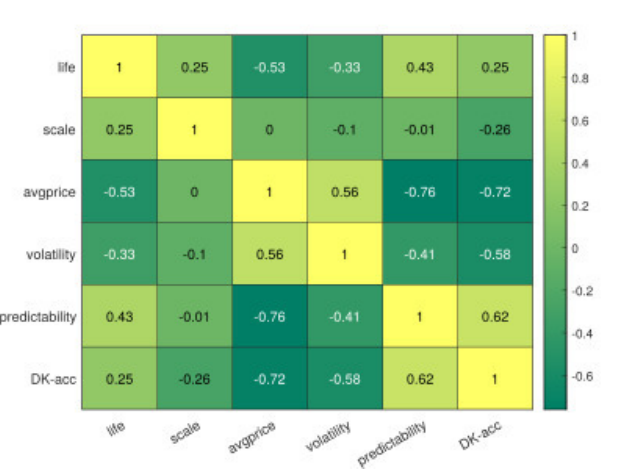}
	\caption{Correlation matrix between different characteristics and prediction accuracy in $T = 0.01$ through DK model.}
	\label{fig:correlation}
\end{figure}
Fig.\ref{fig:correlation} demonstrates the correlation coefficients using Spearman linear regression analysis for the remaining four quantitative features.
Normally, the relationship can be very close with coefficient $[0.7,1.0]$, relatively close $[0.4, 0.7]$, and normal $[0.2,0.4]$, respectively \cite{2011COMPARISON}.
The absolute value correlation coefficients of both $avgprice$ and $volatility$ with ACC are above 0.4, indicating a relatively close correlation, while $life$ and $scale$ are not as relevant.

\begin{figure}
	\centering
	\subfigure[Error bar on accuracy in different price indices]{
	    \label{price-acc}
		\includegraphics[width=0.8\linewidth]{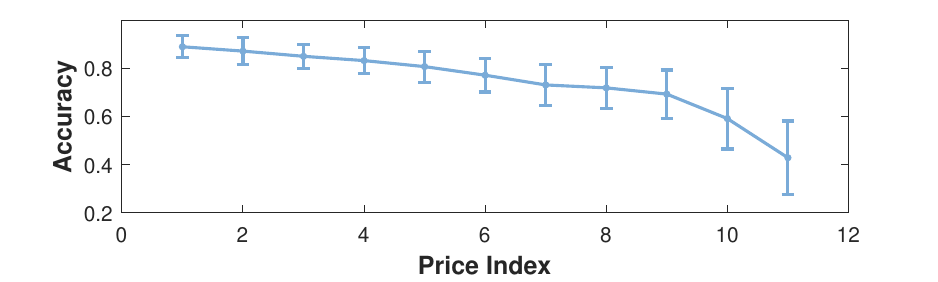}}
	\\
	\subfigure[Error bar on accuracy in different volatility indices]{
	    \label{vol-acc}
		\includegraphics[width=0.8\linewidth]{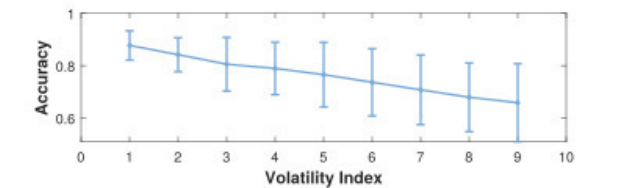} }
	\caption{The error bar on accuracy in different features indices. The horizontal axis represents the index after the classification of features, the vertical axis represents the accuracy of prediction, the points on the I-shaped line segments in the figure means the average prediction accuracy in each index, and the general of the length of the line segment is the standard deviation.}
	\label{feature-acc} 
\end{figure}

Fig.\ref{feature-acc} further presents the details about the impact of $avgprice$ and $volatility$ on the accuracy.
We can find that in the two subplots, the average accuracy of each type of stock gradually decreases as the index increases, while the standard deviation of the prediction accuracy of such stocks gradually increases, which verifies that price and volatility are indeed strongly and negatively correlated with accuracy.

\section{Conclusions}
\label{sec:conclusions}

This paper calculates the predictability of stock prices for tick data of Shanghai and Shenzhen stock markets.
Using MC and DK models for forecasting, we come to the conclusion that there is still a gap between the existing models and the optimal model of theoretical forecasting accuracy, implying that the Chinese stock market is predictable and there is still space to improve current models.
In addition to this, we also discuss the attributes that affect the prediction accuracy and find that the average price and price volatility of the stock have a relatively strong relationship with accuracy.
These results may provide reference value for researchers studying tick data and investors in the Chinese stock market.

There are some aspects in this work that can be improved and added.
We may optimize the models in the future and add more neural network models for comparison. 
Using mathematical models to represent the relationship between prediction accuracy and RMSE is also one of the goals of our future work.
Apart from that, we only considered six profiles when discussing the factors affecting accuracy, there may be other characteristics that we have not tapped into.

\appendix
\section{Details in Discussion}
In this appendix, we provide the details of five features related to the stock data and its belonging enterprises in dataset A. 
We classify the industries of the enterprises to which the stocks belong according to the new China Industry Classification Standard published by Fortune China.
This standard has 23 categories, and the stocks in our data set include only 20 of them, 
which are shown in Table.\ref{industry} for classification details and distribution of industries.
Enterprise scale refers to the number of its employees, and we arbitrarily quantify this number as even as possible, shown in Table.\ref{scale}.
Regions indicate the provinces or equivalent regions where the Enterprise is located, including 22 provinces, 4 municipalities directly under the Central Government and 5 autonomous regions, as well as Hong Kong, Macao and Taiwan regions China (referred to as Outbound) in Table.\ref{region}.
Life indicates how many years the company has been in the stock market by 2021, with the shortest time listed being 1 year and the longest being 31 years.
For non-discrete values like stock price and volatility, we classify them according to their quantitative distribution, as detailed in Table.\ref{price} and Table.\ref{vol}.
All features were normalized before correlation analysis.

\begin{table}
\caption{Industry Classification of the Enterprise}
\centering
\begin{tabular}{|c|c|c|}
Index      & Industry & Number \\
1 & Safety Protection & 15   \\
2 & Office Equipment and Supplies & 39 \\
3 & Electronics and Electrical & 407  \\
4 & Garments and Textiles & 110\\
5 & Environmental Afforestation & 122\\
6 & Mechanical and Electrical Products & 462   \\
7 & Home Supplies & 173   \\
8 & Building and Construction & 282\\
9 & Transportation & 345  \\
10 & Finance & 137   \\
11 & Tourism and Leisure & 35   \\
12 & Agriculture and Farming & 126  \\
13 & Light Industry and Food & 128   \\
14 & Petrochemical Engineering & 376   \\
15 & Water Conservancy and Hydropower & 197   \\
16 & Information Technology & 440   \\
17 & Metallurgy and Mining & 194   \\
18 & Medical and Pharmacy & 357   \\
19 & Professional Services & 151   \\
20 & Comprehensive Industry & 51
\end{tabular}
\label{industry}
\end{table}
\begin{table}
\caption{Quantification of Enterprise Scale}
\centering
\begin{tabular}{|c|c|c|}
Index & Company Scale & Number \\
1 & 0-500 & 513   \\
2 & 500-1,000 & 754 \\
3 & 1,000-1,500 & 507 \\
4 & 1,500-2,000 & 384  \\
5 & 2,000-2,500 & 314   \\
6 & 2,500-3,500 & 356  \\
7 & 3,500-5,000 & 329   \\
8 & 5,000-8,000 & 357   \\
9 & 8,000-15,000 & 297 \\
10 & 15,000+ & 267    \\
\end{tabular}
\label{scale}
\end{table}

\begin{minipage}[c]{0.5\textwidth}
\hspace{-3mm}
    \begin{tabular}{|c|c|c|}
    Index & Region & Number \\
    1 & Anhui & 126   \\
    2 & Beijing & 377 \\
    3 & Fujian & 151  \\
    4 & Gansu & 34\\
    5 & Guangdong & 673\\
    6 & Guangxi & 38   \\
    7 & Guizhou & 30   \\
    8 & Hainan & 33\\
    9 & Hebei & 62  \\
    10 & Henan & 87   \\
    11 & Heilongjiang & 39   \\
    12 & Hubei & 113  \\
    13 & Hunan & 116   \\
    14 & Jilin & 44   \\
    15 & Jiangsu & 479   \\
    16 & Jiangxi & 57   \\
    17 & Outbound & 3   \\
    18 & Liaoning & 75  \\
    19 & Inner Mongolia & 25   \\
    20 & Ningxia & 14   \\
    21 & Qinghai & 10   \\
    22 & Shandong & 228 \\
    23 & Shanxi & 40    \\
    24 & Xi'an & 58 \\
    25 & Shanghai & 338 \\
    26 & Sichuan & 136  \\
    27 & Tianjin & 60   \\
    28 & Xizang & 20    \\
    29 & Xinjiang & 57  \\
    30 & Yunnan & 38    \\
    31 & Zhejiang & 521 \\
    32 & Chongqing & 56
    \end{tabular}
    \captionof{table}{Region List of the Enterprise\label{region}}
\end{minipage}
\hspace{0mm}
\begin{minipage}[c]{0.5\textwidth}
\centering
    \begin{tabular}{|c|c|c|}
    Index & Price(CNY) & Number \\
    1 & 0-3 & 320   \\
    2 & 3-4 & 301 \\
    3 & 4-5 & 304 \\
    4 & 5-6 & 317 \\
    5 & 6-8 & 475   \\
    6 & 8-10 & 374 \\
    7 & 10-13 & 371   \\
    8 & 13-17 & 381   \\
    9 & 17-25 & 400 \\
    10 & 25-50 & 525    \\
    11 & 50+ & 371
    \end{tabular}
    \captionof{table}{Quantification of Stock Price\label{price}}
    \vspace{1cm}
    \begin{tabular}{|c|c|c|}
    Index & Volatility & Number \\
    1 & 0-0.01 & 500   \\
    2 & 0.01-0.016& 491 \\
    3 & 0.016-0.021 & 459 \\
    4 & 0.021-0.027 & 516  \\
    5 & 0.027-0.033 & 480   \\
    6 & 0.033-0.04 & 486  \\
    7 & 0.04-0.05 & 503   \\
    8 & 0.05-0.07 & 449   \\
    9 & 0.07+ & 236 \\
    \end{tabular}
    \captionof{table}{Quantification of Volatility\label{vol}}

\end{minipage}

\newpage
\section*{Acknowledgment}
This work was partially supported by Key Program of Natural Science Foundation of China under Grant(61631018), Huawei Technology Innovative Research.
Yueshan Chen also thanks Miss Fan Zhang for data preprocessing.

\end{document}